\documentclass[a4paper,11pt]{article}
\usepackage[english]{babel}
\usepackage[latin1]{inputenc}   
\usepackage[T1]{fontenc}        

 \usepackage{Vmargin}
 \setpapersize{A4}

\usepackage{amssymb,amsmath,amscd}
\begin{document}

\newcount\cstecnt\cstecnt=0
\def\cste#1{\global\advance\cstecnt by 1
\def\Ht{\mathcal{H}}
\lacste{#1}}
\def\lacste#1{c^{(\the\cstecnt)}_{#1}}

\def\petit#1{\hbox to 0pt{$\DS #1$}}

\def\erfc{\mbox{erfc}}

\def\minx{\min_i X^{(i)}}
\def\minxp{\min_i X'^{(i)}}
\def\maxx{\max_i X^{(i)}}
\def\maxxp{\max_i X'^{(i)}}
\def\argmin{\mbox{arg\hspace{.2em}min}}
\let\PsachantD=\Psachant
\newtheorem{definition}{Définition}
\newtheorem{definitioneng}{Definition}
\renewcommand{\thedefinition}{}
\newtheorem{lemme}{Lemme}
\newtheorem{prop}[lemme]{Proposition}
\newtheorem{theoreme}{Théorème}
\newtheorem{theorem}{Theorem}
\newtheorem{lemma}{Lemma}
\newtheorem{remark}[lemme]{Remark}
\newtheorem{remarque}[lemme]{Remarque}
\newcommand{\di}{{\rm d}}
\def\indic#1{{\rm\bf 1}_{\left\{#1\right\}}}
\def\indice#1{{\rm\bf 1}_{#1}}
\def\preuvenameeng{Proof}
\newtheorem{coroll}[lemme]{Corollary}
\newenvironment{proof}[1][]%
{\noindent{\bf\preuvenameeng{}#1:}\par\nobreak}{\nopagebreak\hspace*{\fill} $\square$\vspace*{1.4ex}\par}

\bibliographystyle{plain}

\title{Asymptotics of a two-dimensional sticky random walk}

\author{Jean B\'erard}

\maketitle
\begin{center}Ex-Laboratoire de
Probabilit\'es, Combinatoire et Statistique 

 Universit\'e Claude
Bernard, Lyon I 

50, avenue Tony Garnier
69366 Lyon Cedex 07

{\tt jean.berard@univ-lyon1.fr}
\end{center}
\abstract{We study the asymptotic behavior of a Markov chain on $\mathbb{Z}^2$ that  corresponds to the two-dimensional marginals of a reinforcement process on $\mathbb{Z}^{\mathbb{N}}$. Three distinct asymptotic regimes are identified, depending on the scaling of the reinforcement parameter $\Delta$ with respect to the number of steps performed by the chain. }

\vspace{1ex}

\vspace{1ex}
\noindent{\sc A.M.S. Subject Classification: 60F05, 60J10, 60K35}.

\vspace{1ex}

\bibliographystyle{plain}

\section{Introduction}

We study the Markov chain $(S_{\Delta}(n))_{n \geq 0} = \left(S^{(1)}_{\Delta}(n), S^{(2)}_{\Delta}(n)   \right)_{n \geq 0}$ on $\mathbb{Z}^2$ defined by the initial condition $\left(S^{(1)}_{\Delta}(0), S^{(2)}_{\Delta}(0)   \right) = (0,0)$ and the following transition kernel $K_{\Delta}[\cdot, \cdot]$ on $\mathbb{Z}^2$:
$$\left\{ \begin{array}{lll} \vspace{2ex} K_{\Delta} \left[ (x,y),(x \pm 1,y \pm 1) \right]= 1/4 \mbox{ for $x \neq y$}, \\ \vspace{2ex}
K_{\Delta} \left[(x,x) , (x+1,x+1) \right] = K_{\Delta} \left[(x,x) , (x-1,x-1) \right] = \frac{u}{4} , \\ \vspace{2ex}

 K_{\Delta} \left[(x,x) , (x+1,x-1) \right] = K_{\Delta} \left[(x,x) , (x+1,x-1) \right] = \frac{2-u}{4},           
   \end{array}  \right. , $$
where $$u = \frac{2+2 \Delta}{2+\Delta},$$
and where $\Delta$ is a fixed non-negative parameter.

The Markov chain $(S_{\Delta}(n))_{n \geq 0}$ defines a sticky random walk on $\mathbb{Z}^2$, in the sense that, when $S^{(1)}_{\Delta}(n) = S^{(2)}_{\Delta}(n)$, the two trajectories are more likely to stick together for the next step than to go apart from each other. Clearly from the definition above, the larger $\Delta$, the stronger the stickiness of the walk. In fact, $(S_{\Delta}(n))_{n \geq 0}$ appears as the two-dimensional marginals of a more complex sticky Markov chain on $\mathbb{Z}^{\mathbb{N}}$, as described in~\cite{Ber}. A similar model, where $\mathbb{Z}$ is replaced by (discretized versions of) the circle $\mathbb{R} / \mathbb{Z}$ has been introduced in~\cite{LeJLem}. From~\cite{Ber}, the behavior of $S_{\Delta}(n)$ for large $n$ and fixed $\Delta$ is degenerate, in the sense that $n^{-1/2} S_{\Delta}(n)$ converges in distribution to a couple of independent standard normal random variables, as would a couple of independent (normalized by $n^{-1/2}$) random walks on $\mathbb{Z}$ do. The stickiness of the walk appears only in the corrections to this first-order behavior. In~\cite{LeJLem}, the $\mathbb{R}/\mathbb{Z}$-version of the model is studied under regimes where $n$ and $\Delta$ go to infinity simultaneously, with $\Delta \sim \alpha n^{1/2}$, and the walk is shown to converge in distribution, after appropriate rescaling, to a non-trivial limiting process, formerly introduced in~\cite{LeJRai}, and known as the sticky flow on the circle. In the present paper, a similar scaling for the two-dimensional marginals of the model on $\mathbb{Z}$ is investigated. In our opinion, the main interest of the present work is the explicit formula obtained for the Fourier transform of the limiting distribution of the walk, since, to our knowledge, it is not possible to derive directly such explicit results from the construction of the sticky flow given in~\cite{LeJLem}. 

Here is our main result.

\begin{theorem}\label{t:cvenloi} 
Assume that the sequence $(\Delta_n)_{n \geq 0}$ goes to infinity as $n$ goes to infinity. The following holds:

\begin{itemize}

\item if $\Delta_n = o(n^{1/2})$, then
$$n^{-1/2} S_{\Delta_n}(n) \stackrel{d}{\longrightarrow} \mathcal{N}(0,1) \otimes \mathcal{N}(0,1),$$

\item if $\Delta_n \sim  \alpha  n^{1/2} $, then  
$$ n^{-1/2} S_{\Delta_n}(n)  \stackrel{d}{\longrightarrow} \mathcal{S}_{\alpha,2}(1),$$
where $\mathcal{S}_{\alpha,2}(1)$ is the probability distribution on $\mathbb{R}^2$ whose Fourier transform
$$ \phi_{\alpha,2}(1)(s,t)=\int_{\mathbb{R}^2} e^{is x_1 + it x_2} d \mathcal{S}_{\alpha,2}(1)(x_1,x_2)$$ 
is given by the following formula
$$ \phi_{\alpha,2}(1)(s,t)= \exp \left( -\frac{s^2+t^2}{2} \right) \left[ 1 - ts \int_0^1  \exp \left(\frac{s^2+t^2}{2} x \right) \ell_{\alpha, s+t}(x)dx  \right], $$
where  $\ell_{\alpha, s+t}$ is defined in Section~\ref{s:ladef} below,

\item if $\Delta_n \gg n^{1/2}$, then  
$$n^{-1/2} S_{\Delta_n}(n) \stackrel{d}{\longrightarrow} (Z,Z),$$
where the random variable $Z$ has the $\mathcal{N}(0,1)$ distribution.

\end{itemize}

\end{theorem}

In the $\Delta_n = o(n^{1/2})$ regime, we conclude that the stickiness of the walk is not strong enough to yield a macroscopic effect on the limiting distribution, since, to the first order, the walk behaves as a couple of independent random walks on $\mathbb{Z}$.
In the $\Delta_n \gg n^{1/2}$ regime, the stickiness is so strong that, to the first order, the walk behaves as a single random walk on $\mathbb{Z}$. The most interesting regime is the $\Delta_n \sim \alpha n^{1/2}$ regime, where the stickiness yields a non-trivial effect on the limiting distribution.  Indeed, in this regime, the distribution of  $n^{-1/2} S_{\Delta_n}(n)$ converges to a non-trivial limiting object, that we identify in~\cite{Ber2} as the distribution of the two-dimensional marginals of the sticky flow on $\mathbb{R}$ at time 1 (yet to be constructed).

The following corollary describes the asymptotic behavior of the covariance of $S^{(1)}_{\Delta_n}(n)$ and $S^{(2)}_{\Delta_n}(n)$. It stems readily from Theorem~\ref{t:cvenloi} above and the (easy) computation of the partial derivative $\frac{\partial^2}{ \partial s \partial t}$ of the Fourier transform of $\mathcal{S}_{\alpha,2}(1)$ at $(t,s)=(0,0)$.

\begin{coroll}
Assume that $\Delta_n \sim  \alpha  n^{1/2} $, then 
$$\lim_{n \to +\infty} n^{-1} E \left[S^{(1)}_{\Delta_n}(n) S^{(2)}_{\Delta_n}(n) \right]  =  \int_{0}^1 \exp \left(\frac{4}{\alpha^2} s \right) \times  \mbox{ \rm erfc} \left(\frac{2}{\alpha} \sqrt{s} \right) ds$$
(the definition of {\rm erfc} is recalled in Section~\ref{s:ladef} below).

\end{coroll}

The structure of the paper is as follows. In Section~\ref{s:ladef}, we define the function  $\ell_{\alpha,w}$ appearing in the expression of the Fourier transform of  $\mathcal{S}_{\alpha,2}(1)$.

In Section~\ref{s:foncgen}, we derive explicit expressions for generating functions related to the walk. In Section~\ref{s:asympt}, an asymptotic analysis of these expressions is performed. In Section~\ref{s:preuveduth}, we use the results of the preceding sections to prove Theorem~\ref{t:cvenloi} above.

\section{Definition of $\ell_{\alpha,w}$}\label{s:ladef}

First, recall that, for all $x \in \mathbb{R}$, 
$$\erfc(x) = 1 - \frac{2}{\sqrt{\pi}} \int_0^x \exp(-y^2) dy.$$
and that this function extends to an entire function on $\mathbb{C}$.

Now, for $w \in \mathbb{R}$ and $\alpha > 0$, such that $\alpha^{-2} - w^2/4 \neq 0$, set 
$$\gamma = \sqrt{ \alpha^{-2} - w^2/4},$$ 
with the convention that when  $\alpha^{-2} - w^2/4 < 0$, $\gamma$ is the square root with positive imaginary part,
and
$$b_1 = -2 \alpha^{-1} + 2 \gamma , \ b_2 =  -2 \alpha^{-1} - 2 \gamma.$$
Define, for all $x \geq 0$, the following function
\begin{eqnarray*}\label{e:defsidiffdezero} && \ell_{\alpha,w} (x) = \\ && \ \  \frac{1}{2 \gamma} \exp \left(-\frac{w^2}{4} x \right) \left[   \frac{b_1}{2} \exp \left( \frac{b_1^2}{4}x \right) \erfc\left(-\frac{b_1}{2} \sqrt{x} \right)  -  \frac{b_2}{2} \exp \left( \frac{b_2^2}{4} x\right) \erfc\left(-\frac{b_2}{2} \sqrt{x} \right)   \right] \ (\ref{e:defsidiffdezero}) .  \end{eqnarray*}

When  $\alpha^{-2} - w^2/4 = 0$, we set
$$\ell_{\alpha,w} (x) =  \frac{1}{2} \exp \left(-\frac{w^2}{4} x \right) \left[   -4 \frac{\sqrt{x}}{ \alpha \sqrt{\pi}} + \left( \frac{4x}{\alpha^2}+2  \right)   \exp\left( \frac{x}{\alpha^2}\right)  \erfc \left( \frac{\sqrt{x}}{\alpha}  \right) \right].$$
This last expression can be obtained as the limit of~(\ref{e:defsidiffdezero}) when $(\alpha,w)$ converges to a limit $(\alpha_*,w_*)$ such that  $\alpha_*^{-2} - w_*^2/4 = 0$.

\section{Generating function computations}\label{s:foncgen}

For all $j \geq 0$, define
$$h_{\Delta}(j,t,n)=\sum_{a \in \mathbb{Z}} e^{ita} P \left(S^{(1)}_{\Delta}(n)=a-j ,S^{(2)}_{\Delta}(n)=a+j \right),$$
and the corresponding generating functions, for $z \in [0,1[$:
$$H_{\Delta}(j,t,z) = \sum_{n=0}^{+\infty} h_{\Delta}(j,t,n) z^n,$$
the above series being well defined since $|h_{\Delta}(j,t,n) | \leq 1$ by definition.

\begin{prop}\label{p:calculexact}
For all $t \in \mathbb{R}$ and $z \in [0,1[$, the following identities hold. 
 $$\frac{1}{H_{\Delta}(0,t,z)}=1-\frac{uz \cos(t)}{2} - (2-u) \left[1-\frac{\cos(t)z}{2} - \sqrt{1-\cos(t) z - \frac{\sin^2(t)}{4} z^2}    \right],$$
$$H_{\Delta}(1,t,z) =H_{\Delta}(0,t,z) \left(\frac{2}{z} - u \cos(t) \right) - \frac{2}{z}  ,$$
and, for all $j \geq 2$,
$$H_{\Delta}(j,t,z) = H_{\Delta}(1,t,z) \left( \frac{2}{z} - \cos(t) - \sqrt{  \left(\frac{2}{z} - \cos(t)\right)^2 -1} \right)^{j-1}.$$
\end{prop}



\begin{proof}[ of Proposition~\ref{p:calculexact}]

A one-step analysis of the Markov chain $(S_\Delta(n))_{n \geq 0}$ yields 
\begin{eqnarray*} && P \left(S^{(1)}_{\Delta}(n+1)=a ,S^{(2)}_{\Delta}(n+1)=a\right)  = \frac{u}{4}  P \left(S^{(1)}_{\Delta}(n)=a-1 ,S^{(2)}_{\Delta}(n)=a-1 \right)  + \\  && \frac{u}{4}  P \left(S^{(1)}_{\Delta}(n)=a+1 ,S^{(2)}_{\Delta}(n)=a+1 \right) +   \frac{1}{4}  P \left(S^{(1)}_{\Delta}(n)=a-1 ,S^{(2)}_{\Delta}(n)=a+1 \right) +  \\ && \frac{1}{4}  P \left(S^{(1)}_{\Delta}(n)=a+1 ,S^{(2)}_{\Delta}(n)=a-1 \right) ,\end{eqnarray*}
whence, taking Fourier transforms on both sides, and using the fact that, for obvious symmetry reasons, $(S^{(1)}_{\Delta}(n),S^{(2)}_{\Delta}(n) )$ and $(S^{(2)}_{\Delta}(n),S^{(1)}_{\Delta}(n) )$ have the same distribution,   
$$ h_{\Delta}(0,t,n+1) =  u\frac{\cos(t)}{2} h_{\Delta}(0,t,n)  + \frac{1}{2} h_{\Delta}(1,t,n).$$
Similarly,
\begin{eqnarray*}    P \left(S^{(1)}_{\Delta}(n+1)=a-1 ,S^{(2)}_{\Delta}(n)=a+1 \right)    = \left(\frac{2-u}{4} \right)   P \left(S^{(1)}_{\Delta}(n)=a ,S^{(2)}_{\Delta}(n)=a \right)     +  \\ \frac{1}{4}  P \left(S^{(1)}_{\Delta}(n)=a-2 ,S^{(2)}_{\Delta}(n)=a+2 \right)  + \frac{1}{4}    P \left(S^{(1)}_{\Delta}(n)=a ,S^{(2)}_{\Delta}(n)=a+2 \right) \\  + \frac{1}{4}   P \left(S^{(1)}_{\Delta}(n)=a-2 ,S^{(2)}_{\Delta}(n)=a \right), 
\end{eqnarray*}
whence
$$  h_{\Delta}(1,t,n+1)  =  \left( \frac{2-u}{4} \right) h_{\Delta}(0,t,n) + \frac{1}{4}  h_{\Delta}(2,t,n) + \frac{\cos(t)}{2}  h_{\Delta}(1,t,n).$$
Finally, for all $j \geq 2$:
\begin{eqnarray*}      P \left(S^{(1)}_{\Delta}(n+1)=a-j ,S^{(2)}_{\Delta}(n+1)=a+j \right)   =     \frac{1}{4} P \left(S^{(1)}_{\Delta}(n)=a-j-1 ,S^{(2)}_{\Delta}(n)=a+j+1 \right) \\ + \frac{1}{4}  P \left(S^{(1)}_{\Delta}(n)=a-j+1 ,S^{(2)}_{\Delta}(n)=a+j-1 \right)   + \frac{1}{4}  P \left(S^{(1)}_{\Delta}(n)=a-j-1 ,S^{(2)}_{\Delta}(n)=a+j-1 \right) \\ +  \frac{1}{4}   P \left(S^{(1)}_{\Delta}(n)=a-j+1 ,S^{(2)}_{\Delta}(n)=a+j+1 \right),
\end{eqnarray*}
whence
$$h_{\Delta}(j,t,n+1)  = \frac{1}{4} h_{\Delta}(j-1,t,n) +  \frac{1}{4} h_{\Delta}(j+1,t,n) + \frac{\cos(t)}{2} h_{\Delta}(j,t,n).$$

Taking generating functions on both sides of the above identities, we obtain that
$$\frac{H_{\Delta}(0,t,z)-1}{z} = u \frac{\cos(t)}{2} H_{\Delta}(0,t,z) + \frac{1}{2} H_{\Delta}(1,t,z),$$
that
$$\frac{H_{\Delta}(1,t,z)}{z} = \frac{2-u}{4} H_{\Delta}(0,t,z) + \frac{1}{4} H_{\Delta}(2,t,z) +\frac{\cos(t)}{2} H_{\Delta}(1,t,z),$$
and, for all $j \geq 2$, that
$$\frac{H_{\Delta}(j,t,z)}{z} = \frac{1}{4} H_{\Delta}(j-1,t,z) +\frac{1}{4} H_{\Delta}(j+1,t,z) +\frac{\cos(t)}{2}H_{\Delta}(j,t,z).$$

For fixed $z$, the sequence  $(H_{\Delta}(j,t,z))_{j \geq 1}$ satisfies an order-two linear induction equation with constant coefficients, whose characteristic equation reads
$$X^2+\left(2 \cos(t) - \frac{4}{z} \right) X +1 = 0. $$
For $0 < z < 1$, the characteristic equation has two distinct real zeros:
$$q_1(z) = \frac{2}{z} - \cos(t) + \sqrt{  \left(\frac{2}{z} - \cos(t)\right)^2 -1}$$ and $$q_2(z) = \frac{2}{z} - \cos(t) - \sqrt{  \left(\frac{2}{z} - \cos(t)\right)^2 -1}.$$
As a consequence, for all $0 < z < 1$, there exist $\alpha(z),\beta(z) \in \mathbb{R}$ such that, for every $j \geq 1$,
\begin{equation}\label{e:reccoeff}H_{\Delta}(j,t,z) = \alpha(z) q_1(z)^{j-1} + \beta(z) q_2(z)^{j-1}.\end{equation}
Observe that $0 < q_2(z) < 1 < q_1(z)$, and that, for fixed $z \in [0,1[$, $|H_{\Delta}(j,t,z)|$ is bounded above by $(1-z)^{-1}$, since $|h_{\Delta}(j,t,n)| \leq 1$ for all $n,j,t$. As a consequence,  one must have $\alpha(z) = 0$, since otherwise, from Identity~(\ref{e:reccoeff}) above, $H_j(z)$ would go to infinity as $j$ goes to infinity. Writing $$H_{\Delta}(1,t,z) = \alpha(z) + \beta(z)$$ and $$H_{\Delta}(2,t,z) = \alpha(z) q_1(z) + \beta(z)q_2(z),$$ we deduce that
$$\beta(z)= H_{\Delta}(1,t,z) ,$$
whence 
$$ H_{\Delta}(j,t,z) = H_{\Delta}(1,t,z)   q_2(z)^{j-1},$$
for all $j \geq 1$. 
Turning back to $H_{\Delta}(0,t,z)$ we can express $H_{\Delta}(2,t,z)$  as a function of $H_{\Delta}(1,t,z)$ in the preceding equations, whence a 2x2 linear system:
$$\frac{H_{\Delta}(0,t,z) (z)-1}{z} = u \frac{\cos(t)}{2} H_{\Delta}(0,t,z) + \frac{1}{2} H_{\Delta}(1,t,z) ,$$
$$\frac{ H_{\Delta}(1,t,z)           }{z} = \frac{2-u}{4} H_{\Delta}(0,t,z)  + \frac{1}{4} q_2(z) H_{\Delta}(1,t,z)  +\frac{\cos(t)}{2}H_{\Delta}(1,t,z).$$
We deduce the expressions of $H_{\Delta}(0,t,z)$ and $H_{\Delta}(1,t,z)$ stated in the Lemma.  

\end{proof}

Introduce the Fourier transform $f_{\Delta}(s,t,n)$  defined by
$$f_{\Delta}(s,t,n) = E\left[\exp \left(is S^{(1)}_{\Delta}(n) + it S^{(2)}_{\Delta}(n)  \right) \right].$$

We now show that  $f_{\Delta}(s,t,n)$ can be expressed in terms of the Fourier transforms $h_{\Delta}(0,s+t,\cdot)$.

\begin{prop}\label{p:pratique} 
For all $s,t,n,\Delta$, the following identity holds:
\begin{eqnarray*} && f_{\Delta}(s,t,n) =  (\cos(t) \cos(s))^n + \\ && \ \ \left[\frac{u}{2}( \cos(t+s)  - \cos(t) \cos(s))  \right] \sum_{k=0}^{n-1} (\cos(t) \cos(s))^k  h_{\Delta}(0,t+s,n-k-1) .\end{eqnarray*}
\end{prop}

\begin{proof}[ of Proposition~\ref{p:pratique}]

A one-step analysis of the Markov chain $(S_{\Delta}(n))_{n \geq 0}$ yields

\begin{eqnarray*} \label{e:unpas} && f_{\Delta}(s,t,n+1) =  E\left[\exp \left(is S^{(1)}_{\Delta}(n+1) + it S^{(2)}_{\Delta}(n+1)  \right) \right] = \\  && \left[ \frac{u}{2} \cos(t+s) + \left(1-\frac{u}{2} \right) \cos(t) \cos(s) \right]   E\left[\exp \left(is S^{(1)}_{\Delta}(n) + it S^{(2)}_{\Delta}(n)  \right) \indic{  S^{(1)}_{\Delta}(n) =  S^{(2)}_{\Delta}(n)  } \right] \\ &&  + \cos(t) \cos(s) E\left[\exp \left(is S^{(1)}_{\Delta}(n) + it S^{(2)}_{\Delta}(n)  \right) \indic{  S^{(1)}_{\Delta}(n) \neq  S^{(2)}_{\Delta}(n)  } \right] \ (\ref{e:unpas}) .\end{eqnarray*}

Observe that, by definition, 
$$h_{\Delta}(0,s+t,n) =   E\left[\exp \left(is S^{(1)}_{\Delta}(n) + it S^{(2)}_{\Delta}(n)  \right) \indic{  S^{(1)}_{\Delta}(n) =  S^{(2)}_{\Delta}(n)  } \right].$$
Using the decomposition $$ \indic{  S^{(1)}_{\Delta}(n) \neq  S^{(2)}_{\Delta}(n)  } = 1 - 
 \indic{  S^{(1)}_{\Delta}(n) =  S^{(2)}_{\Delta}(n)  },$$
Identity~(\ref{e:unpas}) above rewrites
$$f_{\Delta}(s,t,n+1) = \left[\frac{u}{2}( \cos(t+s)  - \cos(t) \cos(s))  \right] h_{\Delta}(0,s+t,n) + \cos(t) \cos(s)f_{\Delta}(s,t,n).$$
Iterating $n$ times the above identity, we obtain the result of the Proposition.

\end{proof}

\section{Asymptotic analysis of generating functions}\label{s:asympt}

In this section, we use the exact formul{\ae}  from the preceding section to prove the following proposition:
\begin{prop}\label{p:laplace}
For all $w \in \mathbb{R}$ and $\theta \geq 0$, the following holds:
\begin{itemize}

\item if $\Delta_n = o(n^{1/2})$, then
\begin{eqnarray*} \lim_{n \to +\infty} n^{-1/2} \Delta_n^{-1} \sum_{k=0}^n h_{\Delta_n}(0,w n^{-1/2},k) \exp\left(\frac{k \theta}{n} \right) = \\ \int_{0}^1 \exp(\theta x) \frac{1}{2} \exp \left(-\frac{w^2}{4}x \right) \frac{2}{\sqrt{\pi x}} dx,\end{eqnarray*}

\item if $\Delta_n \sim  \alpha  n^{1/2} $, then
$$ \lim_{n \to +\infty} n^{-1} \sum_{k=0}^n h_{\Delta_n}(0,w n^{-1/2},k) \exp\left(\frac{k \theta}{n} \right) = \int_{0}^1 \exp(\theta x) \ell_{\alpha,w}(x) dx,$$

\item if $\Delta_n \gg n^{1/2}$, then

$$\lim_{n \to +\infty} n^{-1} \sum_{k=0}^n h_{\Delta_n}(0,w n^{-1/2},k) \exp\left(\frac{k \theta}{n} \right) = \int_{0}^1 \exp(\theta x)  \exp \left(-\frac{w^2}{2}x \right) dx. $$

\end{itemize}

\end{prop}

\begin{proof}[ of Proposition~\ref{p:laplace}]

Define the following family of complex measures on $\mathbb{R}_+$:
$$M_{\Delta_n}(w,n) = n^{-1} \sum_{k=0}^{+\infty}   h_{\Delta_n}(0,w n^{-1/2},k) \delta_{\frac{k}{n}}.$$
The Laplace transform of $M_{\Delta_n}(w,n)$ is defined, for all $\lambda > 0$, by
$$\mathcal{L} \left[M_{\Delta_n}(w,n) \right] ( \lambda ) = \int_{0}^{+\infty} e^{-\lambda x} d M_{\Delta_n}(w,n)(x),$$
so that:
$$   \mathcal{L} \left[ M_{\Delta_n}(w,n) \right] ( \lambda )     =n^{-1} H_{\Delta_n} \left(0,wn^{-1/2} , \exp{(-\lambda n^{-1})} \right).$$

Moreover, we have
\begin{equation}\label{e:superutile}n^{-1} \sum_{k=0}^n h_{\Delta_n}(0,w n^{-1/2},k) \exp\left(\frac{k \theta}{n} \right) = \int_{0}^1 \exp( \theta x) d M_{\Delta_n}(w,n).\end{equation}

\subsubsection*{The $\Delta_n \sim \alpha n^{1/2}$ regime}

An asymptotic analysis of the explicit expression of Proposition~\ref{p:calculexact} proves that 
$$ \lim_{n \to +\infty}   \mathcal{L} \left[ M_{\Delta_n}(w,n) \right] ( \lambda ) =  \frac{1}{w^2/2 + \lambda + \alpha^{-1} \sqrt{ 4 \lambda + w^2}}.$$

We identify the above expression as the Laplace transform at $\lambda$ of the function $\ell_{\alpha,w}(\cdot)$, so we deduce the result of the proposition from standard results on Laplace transforms (see e.g.~\cite{Fel} chap.~13) and from Identity~(\ref{e:superutile}) above.









\subsubsection*{The $\Delta_n \gg n^{1/2}$ regime}

It is easily checked from the explicit expression of Proposition~\ref{p:calculexact} that 
$$ \lim_{n \to +\infty}   \mathcal{L} \left[ M_{\Delta_n}(w,n) \right] ( \lambda ) =  \frac{1}{w^2/2 + \lambda} .$$

We identify the above expression as the Laplace transform at $\lambda$ of the function $\exp \left( -\frac{w^2}{2} x \right)$, and the same standard arguments on Laplace transforms as above apply.

\subsubsection*{The $\Delta_n = o(n^{1/2})$ regime }

Again, an asymptotic analysis yields
$$ \lim_{n \to +\infty}  \mathcal{L} \left[ n^{-1/2} \Delta_n^{-1}  M_{\Delta_n}(w,n) \right] ( \lambda ) = \frac{1}{\sqrt{ 4 \lambda + w^2}},$$whence the result of the proposition follows, identifying this last expression as the Laplace transform at $\lambda$ of the function
$$ \frac{1}{2} \exp \left( -\frac{w^2}{4}x \right) \frac{2}{\sqrt{\pi x}}.$$

\end{proof}

\section{Proof of Theorem~\ref{t:cvenloi}}\label{s:preuveduth}

\subsubsection*{The $\Delta_n \sim \alpha n^{1/2}$ regime}

We start from the identity given in Proposition~\ref{p:pratique} above:

\begin{eqnarray*}\label{e:pratique} && f_{\Delta}(s,t,n) =(\cos(t) \cos(s))^n + \\ && \ \  \left[\frac{u}{2}( \cos(t+s)  - \cos(t) \cos(s))  \right] \sum_{k=0}^{n-1} (\cos(t) \cos(s))^k  h_{\Delta}(0,t+s,n-k-1) \ (\ref{e:pratique}).\end{eqnarray*}

For all $0 \leq k \leq n$,  
$$(\cos(t n^{-1/2}) \cos(s n^{-1/2}))^k =  \exp\left[ (n-k) \frac{t^2+s^2}{2n} - \frac{t^2+s^2}{2} \right] \left( 1 + O\left(\frac{k}{n^{3/2}} \right) \right).$$

On the other hand, 
$$\frac{u_n}{2}( \cos(tn^{-1/2}+sn^{-1/2}) - \cos(tn^{-1/2}) \cos(sn^{-1/2}))  =- \frac{ts}{n} + o(n^{-3/2}),$$
where $u_n =  \frac{2+2 \Delta_n}{2+\Delta_n}$.

As a consequence, Equation~(\ref{e:pratique}) above entails that 
\begin{eqnarray*} && f_{\Delta_n}(s n^{-1/2}, t n^{-1/2} ,n ) =  \exp\left(- \frac{s^2+t^2}{2} \right) \\ && - \frac{st}{n} \times \exp \left(-\frac{t^2+s^2}{2} \right) \sum_{k=0}^{n-1} h_{\Delta_n}(0,(t+s)n^{-1/2},n-k-1)  \exp\left[(n-k) \frac{t^2+s^2}{2n} \right] + o(1),\end{eqnarray*}
using the bound $| h_{\Delta}|  \leq 1$ to handle the error terms.

Applying Proposition~\ref{p:laplace} above with  $w=  (s+t)$ and  $\theta=\frac{s^2+t^2}{2}$, we see that 
\begin{eqnarray*}\label{e:labellerouge} && \lim_{n \to + \infty} f_{\Delta_n}(s n^{-1/2}, t n^{-1/2} ,n ) = \\ && \ \   \exp \left( -\frac{s^2+t^2}{2} \right) \left[ 1 - ts \int_0^1  \exp \left(\frac{s^2+t^2}{2} x \right) \ell_{\alpha, s+t}(x)dx  \right] \ (\ref{e:labellerouge}).\end{eqnarray*}
Observe that the sequence of random variables $(n^{-1/2} S_{\Delta_n}(n))_{n \geq 0}$ is tight in distribution, since, for every $\Delta \geq 0$, $S^{(1)}_{\Delta}(n)$ and  $S^{(2)}_{\Delta}(n)$ are both distributed as the position at time $n$ of a simple symmetric random walk on $\mathbb{Z}$. 
By the standard result on the convergence of Fourier transforms (see e.g.~\cite{Shi} chap.~3) and Identity~\ref{e:labellerouge} above, we deduce the result of the theorem in the $\Delta_n \sim \alpha n^{1/2}$ regime, that is, we deduce that $n^{-1/2} S_{\Delta_n}(n)$ converges to a  limiting distribution $\mathcal{S}_{\alpha,2}(1)$ on $\mathbb{R}^2$, and that the Fourier transform of  $\mathcal{S}_{\alpha,2}(1)$ is indeed given by the formula of the theorem.

\subsubsection*{The $\Delta_n \gg n^{1/2}$ regime}
 
The proof is similar, the only difference being that the limiting Fourier transform we obtain from Proposition~\ref{p:laplace} is not the same:
$$ \lim_{n \to + \infty} f_{\Delta_n}(s n^{-1/2}, t n^{-1/2} ,n ) =  \exp \left( -\frac{(s+t)^2}{2} \right).$$

\subsubsection*{The $\Delta_n = o(n^{1/2})$ regime}

The proof is similar. From Proposition~\ref{p:laplace},
$$ \lim_{n \to + \infty} n^{-1} \sum_{k=0}^n (\cos(t) \cos(s))^k  h_{\Delta}(0,t+s,n-k) = 0,$$ so the only term yielding a non-negligible contribution in Identity~\ref{e:pratique} above is $$\left(\cos(tn^{-1/2}) \cos(s n^{-1/2})\right)^n,$$ 
whence
$$ \lim_{n \to + \infty} f_{\Delta_n}(s n^{-1/2}, t n^{-1/2} ,n ) =  \exp \left( -\frac{s^2+t^2}{2} \right).$$

This ends the proof of the theorem. \nopagebreak\hspace*{\fill} $\square$\vspace*{1.4ex}\par

\bibliography{sticky2d.bbl}

\end{document}